\documentclass{amsart}
\usepackage{amsfonts}
\usepackage{cite}

\newtheorem{theorem}{Theorem}
\theoremstyle{plain}

\newtheorem{corollary}{Corollary}

\newtheorem{proposition}{Proposition}
\newtheorem{remark}{Remark}

\numberwithin{equation}{section}

\begin{document}
\title{On generalized quasi-Einstein GRW space-times}
\author{Uday Chand De}
\address[U. C. De]{ Department of Pure Mathematics, University of Calcutta\\
35, Ballygaunge Circular Road\\
Kolkata 700019, West Bengal, India,}
\email{uc\_de@yahoo.com}
\author{Sameh Shenawy}
\address[S. Shenawy]{Basic Science Department, Modern Academy for
Engineering and Technology, Maadi, Egypt,}
\email{drshenawy@mail.com}
\subjclass[2000]{Primary 53C25; Secondary 83F05}
\keywords{Quasi-Einstein manifolds, Perfect fluid space-times, Generalized
Robertson-Walker space-times, Einstein manifolds}

\begin{abstract}
Recently, it is proven that generalized Robertson-Walker space-times in all
orthogonal subspaces of Gray's decomposition but one(unrestricted) are
perfect fluid space-times. GRW space-times in the unrestricted subspace are
identified by having constant scalar curvature. Generalized quasi-Einstein
GRW space-times have a constant scalar curvature. It is shown that
generalized quasi-Einstein GRW space-times reduce to Einstein space-times or
perfect fluid space-times.
\end{abstract}

\maketitle

\section{Introduction}

The warped product $M=I\times _{f}M^{\ast }$ of an open connected interval $%
\left( I,-\mathrm{d}t^{2}\right) $ of $%
\mathbb{R}
$ and a Riemannian manifold $M^{\ast }$ with warping function $%
f:I\rightarrow 
\mathbb{R}
^{+}$ is called a generalized Robertson-Walker space-time(or GRW space-times)%
\cite{Mantica:2017,Sanches:1998}. This family of Lorentzian space-times
broadly extends the classical Robertson-Walker space-times, Friedmann
cosmological models, Einstein-de Sitter space-times and many others\cite%
{Chen:2014,Sanches:1998}. The classical Robertson-Walker spacetime is
regarded as cosmological models since it is spatially homogeneous and
spatially isotropic whereas GRW space-times serve as inhomogeneous extension
of Robertson-Walker space-times that admit an isotropic radiation\cite%
{Chen:2014}(see also \cite{Ehlers:1968,Sanches:1998}). A Lorentzian manifold
is called a perfect fluid space-time if the Ricci tensor $\mathrm{Ric}$
takes the form%
\begin{equation*}
\mathrm{Ric}\left( X,Y\right) =ag\left( X,Y\right) +bA\left( X\right)
A\left( Y\right) 
\end{equation*}%
where $a,b$ are scalars and $A$ is a $1-$form metrically equivalent to a
unit time-like vector field\cite{Mantica:2019,Mantica:2019B}. Perfect fluid
space-times in the language of differential geometry are called
quasi-Einstein spaces where $A$ is metrically equivalent to a unit
space-like vector field. Recently, in \cite{Mantica:2019B}, it is proven
that a perfect fluid space-time with divergence-free conformal curvature
tensor is a GRW space-time with Einstein fibers given that the scalar
curvature is constant. Many sufficient conditions on perfect fluid
space-times to be a GRW space-time are derived.

Gray presented an invariant orthogonal decomposition of the covariant
derivative of the Ricci tensor\cite{Gray:1978}(see also \cite{Mantica:2012}%
). Recently, Carlo Mantica et al proved that the Ricci tensor of a
generalized Robertson-Walker space-time in all classes of Gray's
decomposition but $\mathcal{A\oplus B}$ is either Einstein or takes the form
of a perfect fluid whereas $\mathcal{A\oplus B}$ is not restricted\cite%
{Mantica:2019}. The class $\mathcal{A\oplus B}$ is characterized by $\nabla
R=0$ i.e. the scalar curvature is constant. Now, the following question
arises.

Does the Ricci tensor of all GRW space-times in $\mathcal{A\oplus B}$ reduce
to be Einstein or take the form of a perfect fluid?

In this work, we get a partial positive answer. A (pseudo-)Riemannian
manifold $\left( M,g\right) $ is called a generalized quasi-Einstein
manifold if its Ricci curvature satisfies%
\begin{equation}
\mathrm{Ric}\left( X,Y\right) =\alpha g\left( X,Y\right) +\beta A\left(
X\right) A\left( Y\right) +\gamma \left[ A\left( X\right) B\left( Y\right)
+A\left( Y\right) B\left( X\right) \right] ,  \label{Q0}
\end{equation}%
where $\alpha ,\beta $ and $\gamma $ are non-zero constants, $A$ and $B$ are 
$1-$forms corresponding to two orthonormal vector field\cite%
{Chaki:2001,De:2014,Guler:2015,Guler:2016A,Guler:2016,Mallick:2016}. If $%
\gamma =0$, then $M$ reduces to a quasi-Einstein manifold. It is clear that
generalized quasi-Einstein space-times are generally imperfect fluid
space-times with constant scalar curvature $R=n\alpha +\beta $. However, we
prove that generalized quasi-Einstein GRW space-times are either Einstein or
perfect fluid space-times.

\begin{remark}
It is noted that any vector field orthogonal to a time-like vector field is
space-like. Thus the generators couldn't be time-like. Now, one may assume
that one of the generators is time-like and the other is space-like.
Generally, the results of this article still hold in this case with minor
changes.
\end{remark}

\section{Notes on generalized quasi-Einstein manifolds}

Let $M$ be a generalized quasi-Einstein (pseudo-)Riemannian manifold i.e.%
\begin{equation}
R_{ij}=\alpha g_{ij}+\beta A_{i}A_{j}+\gamma \left(
A_{i}B_{j}+A_{j}B_{i}\right) ,  \label{Q1}
\end{equation}%
where $A_{i}A^{i}=B_{j}B^{j}=1$ and $A_{i}B^{i}=0$. The trace of this
equations gives%
\begin{equation*}
R=n\alpha +\beta .
\end{equation*}%
It is noted that%
\begin{eqnarray*}
A^{i}A^{j}R_{ij} &=&\alpha +\beta , \\
A^{i}B^{j}R_{ij} &=&\gamma , \\
B^{i}B^{j}R_{ij} &=&\alpha ,
\end{eqnarray*}%
and consequently $\beta =\left( A^{i}A^{j}-B^{i}B^{j}\right) R_{ij}$. Hence
we can state the following result.

\begin{proposition}
Let $M$ be a generalized quasi-Einstein manifold with generators $A$ and $B$%
. Then the scalar curvature is%
\begin{equation*}
R=\left( n-1\right) B^{i}B^{j}R_{ij}+A^{i}A^{j}R_{ij}.
\end{equation*}
\end{proposition}

Now assume that $A$ is an eigenvector of the Ricci tensor with eigenvalue $%
\xi $ i.e. $A^{i}R_{ij}=\xi A_{j}$. A contraction of Equation (\ref{Q1})
with $A^{i}$ yields%
\begin{equation*}
A^{i}R_{ij}=\alpha A^{i}g_{ij}+\beta A^{i}A_{i}A_{j}+\gamma \left(
A^{i}A_{i}B_{j}+A^{i}A_{j}B_{i}\right)
\end{equation*}%
implies $\left( \xi -\alpha -\beta \right) A_{j}=\gamma B_{j}$. Thus $\gamma
=0$ and $\xi =\alpha +\beta $. If $B$ is an eigenvector of the Ricci tensor
with eigenvalue $\phi $, then%
\begin{equation*}
B^{i}R_{ij}=\alpha B^{i}g_{ij}+\beta B^{i}A_{i}A_{j}+\gamma \left(
B^{i}A_{i}B_{j}+B^{i}A_{j}B_{i}\right)
\end{equation*}%
infers%
\begin{equation*}
\left( \phi -\alpha \right) B_{j}=\gamma A_{j}.
\end{equation*}%
Consequently, $\phi =\alpha $ and $\gamma =0$. Conversely, assume that $M$
is a quasi-Einstein manifold with generator $A$. Then%
\begin{equation*}
R_{ij}=\alpha g_{ij}+\beta A_{i}A_{j}
\end{equation*}%
yields%
\begin{equation*}
A^{i}R_{ij}=\alpha A^{i}g_{ij}+\beta A^{i}A_{i}A_{j}=\left( \alpha +\beta
\right) A_{j}.
\end{equation*}%
Also,%
\begin{equation*}
B^{i}R_{ij}=\alpha B^{i}g_{ij}+\beta B^{i}A_{i}A_{j}=\alpha A_{j}.
\end{equation*}%
This leads to the following.

\begin{theorem}
Let $M$ be a generalized quasi-Einstein manifold. Then, $M$ reduces to a
quasi-Einstein manifold if and only if one of the generators is an
eigenvector of the Ricci tensor.
\end{theorem}

The covariant derivative of the Ricci tensor of a generalized quasi-Einstein
manifold is given by%
\begin{eqnarray*}
\nabla _{k}R_{ij} &=&\beta \left( \nabla _{k}A_{i}\right) A_{j}+\beta
A_{i}\left( \nabla _{k}A_{j}\right) \\
&&+\gamma \left[ \left( \nabla _{k}A_{i}\right) B_{j}+A_{i}\left( \nabla
_{k}B_{j}\right) +\left( \nabla _{k}A_{j}\right) B_{i}+A_{j}\left( \nabla
_{k}B_{i}\right) \right]
\end{eqnarray*}%
and hence%
\begin{eqnarray*}
\nabla _{i}R_{kj} &=&\beta \left( \nabla _{i}A_{k}\right) A_{j}+\beta
A_{k}\left( \nabla _{i}A_{j}\right) \\
&&+\gamma \left[ \left( \nabla _{i}A_{k}\right) B_{j}+A_{k}\left( \nabla
_{i}B_{j}\right) +\left( \nabla _{i}A_{j}\right) B_{k}+A_{j}\left( \nabla
_{i}B_{k}\right) \right]
\end{eqnarray*}%
Thus, the Codazzi deviation tensor $\mathcal{D}$ is%
\begin{eqnarray*}
\mathcal{D}_{kij} &=&\nabla _{k}R_{ij}-\nabla _{i}R_{kj} \\
&=&\beta \left( \nabla _{k}A_{i}\right) A_{j}+\beta A_{i}\left( \nabla
_{k}A_{j}\right) -\beta \left( \nabla _{i}A_{k}\right) A_{j}-\beta
A_{k}\left( \nabla _{i}A_{j}\right) \\
&&+\gamma \left[ \left( \nabla _{k}A_{i}\right) B_{j}+A_{i}\left( \nabla
_{k}B_{j}\right) +\left( \nabla _{k}A_{j}\right) B_{i}+A_{j}\left( \nabla
_{k}B_{i}\right) \right] \\
&&-\gamma \left[ \left( \nabla _{i}A_{k}\right) B_{j}+A_{k}\left( \nabla
_{i}B_{j}\right) +\left( \nabla _{i}A_{j}\right) B_{k}+A_{j}\left( \nabla
_{i}B_{k}\right) \right]
\end{eqnarray*}

Now, we have the following cases%
\begin{eqnarray*}
A^{j}\mathcal{D}_{kij} &=&\beta \left( \nabla _{k}A_{i}-\nabla
_{i}A_{k}\right) +\gamma \left( \nabla _{k}B_{i}-\gamma \nabla
_{i}B_{k}\right) \\
B^{j}\mathcal{D}_{kij} &=&\nabla _{k}A_{i}-\nabla _{i}A_{k}
\end{eqnarray*}%
Thus, for a Codazzi Ricci tensor, the generators are both closed. In this
case,%
\begin{eqnarray*}
0 &=&\mathcal{D}_{kij} \\
&=&\beta A_{i}\left( \nabla _{k}A_{j}\right) -\beta A_{k}\left( \nabla
_{i}A_{j}\right) +\gamma A_{i}\left( \nabla _{k}B_{j}\right) \\
&&+\gamma \left( \nabla _{k}A_{j}\right) B_{i}-\gamma A_{k}\left( \nabla
_{i}B_{j}\right) -\gamma \left( \nabla _{i}A_{j}\right) B_{k} \\
&=&\nabla _{j}\left[ \beta A_{i}A_{k}+\gamma A_{i}B_{k}+\gamma A_{k}B_{i}%
\right] \\
&&-2\beta A_{k}\left( \nabla _{j}A_{i}\right) -2\gamma \left( \nabla
_{j}A_{i}\right) B_{k}-2\gamma A_{k}\left( \nabla _{j}B_{i}\right) \\
&=&\nabla _{j}R_{ik}-2\left( \beta A_{k}+\gamma B_{k}\right) \nabla
_{j}A_{i}-2\gamma A_{k}\nabla _{j}B_{i}.
\end{eqnarray*}%
Thus we have the following.

\begin{proposition}
Let $M$ be a generalized quasi-Einstein manifold. Assume that $M$ is
Einstein-like of class $\mathcal{B}$(i.e. the Ricci tensor is a Codazzi
tensor). Then $A$ and $B$ are closed. Moreover,%
\begin{equation*}
\nabla _{j}R_{ik}=2\left( \beta A_{k}+\gamma B_{k}\right) \nabla
_{j}A_{i}+2\gamma A_{k}\nabla _{j}B_{i}.
\end{equation*}
\end{proposition}

A contraction of $\mathcal{D}_{kij}$ by $g^{ij}$ and then by the generators $%
A^{k}$ and $B^{k}$ infers%
\begin{eqnarray*}
0 &=&\left( \beta A^{i}+\gamma B^{i}\right) \nabla _{i}A_{k}+\left( \beta
A_{k}+\gamma B_{k}\right) \left( \nabla _{i}A^{i}\right) +\gamma A_{k}\left(
\nabla _{i}B^{i}\right) +\gamma A^{i}\nabla _{i}B_{k}, \\
0 &=&\beta \left( \nabla _{i}A^{i}\right) +\gamma \nabla _{i}B^{i}, \\
0 &=&\gamma \nabla _{i}A^{i}.
\end{eqnarray*}%
Thus $\nabla _{i}A^{i}=\nabla _{i}B^{i}=0$.

Assume that $M$ is Einstein-like of class $\mathcal{P}$( that is, the Ricci
tensor is a parallel, $\nabla _{k}R_{ij}=0$). Then,%
\begin{eqnarray*}
0 &=&\beta \left( \nabla _{i}A_{k}\right) A_{j}+\beta A_{k}\left( \nabla
_{i}A_{j}\right) +\gamma \left( \nabla _{i}A_{k}\right) B_{j} \\
&&+\gamma A_{k}\left( \nabla _{i}B_{j}\right) +\gamma \left( \nabla
_{i}A_{j}\right) B_{k}+\gamma A_{j}\left( \nabla _{i}B_{k}\right)
\end{eqnarray*}

Contractions by $A^{k}$ and $B^{k}$ imply%
\begin{eqnarray*}
0 &=&\beta \nabla _{i}A_{j}+\gamma \nabla _{i}B_{j}, \\
0 &=&\gamma \nabla _{i}A_{j}.
\end{eqnarray*}

Assume that $A$ is not parallel, then $\beta =\gamma =0$. Thus we conclude.

\begin{theorem}
Let $M$ be a Ricci-symmetric generalized quasi-Einstein manifold. Then, $M$
is Einstein if the generator $A$ is not covariantly constant.
\end{theorem}

\section{Generalized quasi-Einstein GRW space-times}

A Lorentzian manifold $M$ is a GRW space-time if and only if $M$ has a unit
time-like vector field $u_{i}$ such that 
\begin{equation*}
\nabla _{k}u_{j}=\varphi \left( g_{kj}+u_{k}u_{j}\right) ,
\end{equation*}%
which is also an eigenvector of the Ricci tensor i.e. $R_{ij}u^{i}=\xi u_{j}$
for some scalar functions $\varphi $ and $\xi $\cite%
{Mantica:2016,Mantica:2017,Mantica:2019}. We say that $u$ is a nontrivial
torse-forming vector field if $\varphi \neq 0$. This characterization is an
alternative of Chen's theorem in \cite{Chen:2014}. If $M$ is a generalized
quasi-Einstein manifold, then%
\begin{equation*}
R_{ij}=\alpha g_{ij}+\beta A_{i}A_{j}+\gamma \left(
A_{i}B_{j}+A_{j}B_{i}\right) .
\end{equation*}%
A contraction by $u^{i}$ yields%
\begin{equation*}
u^{i}R_{ij}=\alpha u^{i}g_{ij}+\beta u^{i}A_{i}A_{j}+\gamma \left(
u^{i}A_{i}B_{j}+A_{j}u^{i}B_{i}\right) 
\end{equation*}%
which implies%
\begin{equation*}
\left( \xi -\alpha \right) u_{j}=\beta \left( u^{i}A_{i}\right) A_{j}+\gamma
\left( u^{i}A_{i}\right) B_{j}+\gamma A_{j}\left( u^{i}B_{i}\right) 
\end{equation*}%
and hence%
\begin{equation*}
\left( \xi -\alpha \right) u_{j}=\left[ \beta \left( u^{i}A_{i}\right)
+\gamma \left( u^{i}B_{i}\right) \right] A_{j}+\gamma \left(
u^{i}A_{i}\right) B_{j}
\end{equation*}%
Two different contractions by the generators give%
\begin{equation*}
\left( \xi -\alpha -\beta \right) \left( u^{i}A_{i}\right) -\gamma \left(
u^{i}B_{i}\right) =0
\end{equation*}%
and%
\begin{equation*}
\left( \xi -\alpha \right) \left( u^{i}B_{i}\right) -\gamma \left(
u^{i}A_{i}\right) =0.
\end{equation*}%
Thus%
\begin{eqnarray*}
\left( \xi -\alpha \right) u_{j} &=&\left[ \beta \left( u^{i}A_{i}\right)
+\left( \xi -\alpha -\beta \right) \left( u^{i}A_{i}\right) \right]
A_{j}+\left( \xi -\alpha \right) \left( u^{i}B_{i}\right) B_{j} \\
&=&\left( \xi -\alpha \right) \left( u^{i}A_{i}\right) A_{j}+\left( \xi
-\alpha \right) \left( u^{i}B_{i}\right) B_{j}
\end{eqnarray*}%
and hence%
\begin{equation}
\left( \xi -\alpha \right) \left[ u_{j}-\left( u^{i}A_{i}\right)
A_{j}-\left( u^{i}B_{i}\right) B_{j}\right] =0.
\end{equation}%
It is clear that $u_{j}$ is not a linear combination of $A_{j}$ and $B_{j}$
only since $u^{i}$ is time-like whereas $A^{i}$ and $B^{i}$ are orthonormal
space-like fields so $\xi =\alpha $. Therefore,%
\begin{eqnarray}
\beta \left( u^{i}A_{i}\right) +\gamma \left( u^{i}B_{i}\right)  &=&0
\label{Q21} \\
\gamma \left( u^{i}A_{i}\right)  &=&0  \label{Q22}
\end{eqnarray}%
It is noted that $\gamma =\beta =0$ if $\left( u^{i}A_{i}\right) $ is not
zero. Suppose that $\left( u^{i}B_{i}\right) $ does not vanish. Then
Equation (\ref{Q22}) implies that either $\gamma =0$ or $\left(
u^{i}A_{i}\right) =0$. The later case with Equation (\ref{Q21}) yield $%
\gamma =0$ i.e. $M$ is quasi-Einstein if $\left( u^{i}B_{i}\right) \neq 0$.
Now, assume that $u^{i}$ is orthogonal to both the generators i.e. $\left(
u^{i}A_{i}\right) =\left( u^{i}B_{i}\right) =0$. The Ricci tensor of a GQE
manifold is%
\begin{equation*}
R_{ij}=\xi g_{ij}+\beta A_{i}A_{j}+\gamma \left(
A_{i}B_{j}+A_{j}B_{i}\right) 
\end{equation*}%
and so%
\begin{equation*}
\nabla _{k}R_{ij}=\beta A_{j}\nabla _{k}A_{i}+\beta A_{i}\nabla
_{k}A_{j}+\gamma \left( B_{j}\nabla _{k}A_{i}+A_{i}\nabla
_{k}B_{j}+A_{j}\nabla _{k}B_{i}+B_{i}\nabla _{k}A_{j}\right) 
\end{equation*}%
A contraction by $u^{i}$ implies%
\begin{eqnarray*}
u^{i}\nabla _{k}R_{ij} &=&0 \\
\nabla _{k}\left( u^{i}R_{ij}\right) -R_{ij}\nabla _{k}u^{i} &=&0
\end{eqnarray*}%
It is noted that $u^{i}$ is an eigenvector of the Ricci tensor ( i.e. $%
u^{i}R_{ij}=\xi u_{j}$) and $\nabla _{k}u^{i}=\varphi \left( \delta
_{k}^{i}+u_{k}u^{i}\right) $. Thus%
\begin{eqnarray*}
\nabla _{k}\left( \xi u_{j}\right) -R_{ij}\varphi \left( \delta
_{k}^{i}+u_{k}u^{i}\right)  &=&0 \\
\xi \nabla _{k}u_{j}-\varphi \delta _{k}^{i}R_{ij}-\varphi u_{k}u^{i}R_{ij}
&=&0 \\
\xi \varphi \left( g_{kj}+u_{k}u_{j}\right) -\varphi R_{kj}-\varphi \xi
u_{k}u_{j} &=&0 \\
\varphi \left( \xi g_{kj}-R_{kj}\right)  &=&0
\end{eqnarray*}%
So $M$ is Einstein if $u$ is a nontrivial torse-forming vector field.

\begin{theorem}
Let $M$ be a generalized quasi-Einstein GRW space-time. Then $%
u^{i}R_{ij}=\alpha u_{j}$ i.e. $\alpha $ is the eigenvalue of the
eigenvector $u^{i}$ and

\begin{enumerate}
\item $M$ reduces to be Einstein space-time if $u^{i}$ is orthogonal to both
the generators provided $\varphi \neq 0$.

\item $M$ reduces to be Einstein space-time if $u^{i}$ is not orthogonal to
first generator.

\item $M$ reduces to be perfect fluid space-time if $u^{i}$ is not
orthogonal to the second generator.
\end{enumerate}
\end{theorem}

\begin{corollary}
Let $M$ be a generalized quasi-Einstein Lorentzian manifold admitting a unit
time-like non-trivial torse-forming vector field. Then $M$ reduces to an
Einstein GRW space-time or a perfect fluid GRW space-time.
\end{corollary}

\end{document}